\documentclass[12pt]{article}
\font\smallit=cmti10

\usepackage{amsmath}
\usepackage{amssymb}
\usepackage{graphicx}
\newcommand{\norm}[1]{\left\vert#1\right\vert}
\newcommand{\set}[1]{\left\{#1\right\}}

\newcommand{\posint}{{\mathbb Z}^+}
\newcommand{\To}{\mapsto}
\usepackage{mathrsfs}

\usepackage{tikz}
\usetikzlibrary{arrows,snakes,backgrounds}
\begin{document}

\begin{center}
{\bf Counting hexagonal lattice animals confined to a strip}
\vskip 20pt
{\bf Moa Apagodu
\footnote{This article is accompanied by two Maple packages, {\tt HexANIMALS}
and {\tt HexaFreeANIMALS}, that can be downloaded from
{\tt http://www.people.vcu.edu/$\sim$mapagodu/}}}\\
{\smallit Department of Mathematics, Virginia Commonwealth University, Richmond, VA 08854, USA}\\
{\tt mapagodu@vcu.edu}\\
\vskip 10pt
\bf Stirling Chow \\
{\smallit Department of Computer Science, University of Victoria, Victoria, BC V8W3P6, Canada}\\
{\tt schow@cs.uvic.ca}\\
\end{center}

\vskip 10pt

\centerline{\bf Abstract}

We describe a bijection between hexagonal lattice animals and a special type of square lattice animals. Using this bijection we adopt Maple packages that automatically generates generating functions (and series expansions) for fixed square lattice lattice animals to that of fixed hexagonal animals on the two-dimensional hexagonal lattice confined to a strip $0\leq y\leq k$, for arbitrary $k$.\\

Polyominoes, Hexagonal Lattice, Square Lattice, Enumeration.
\thispagestyle{empty}
\baselineskip=15pt
\vskip 30pt

\section*{\normalsize 1. Background}

\textbf{Required reading}: Symbol-Crunching with the Transfer-Matrix
Method in Order to Count Skinny Physical Creatures [Z1].\\

In [Z1], Zeilberger used finite transfer matrix method and developed two Maple packages {\tt ANIMALS} and {\tt
FreeANIMALS} to count square lattice animals confined to a strip. Here we define a
bijection between hexagonal lattice animals and a special class of square lattice animals and adopt the
packages {\tt HexANIMALS} and {\tt HexaFreeANIMALS} accompanying [Z1] to enumerate hexagonal lattice animals.

{\bf Definition 1.1} A hexagonal lattice animal (hexagonal polyominoes hence forth) on the hexagonal lattice is an edge-connected set of lattice cells on the hexagonal lattice.  Two animals are equivalent if they are translations of each other.\\

Figure 1 shows the three hexagonal polyominoes comprised of two hexagons.

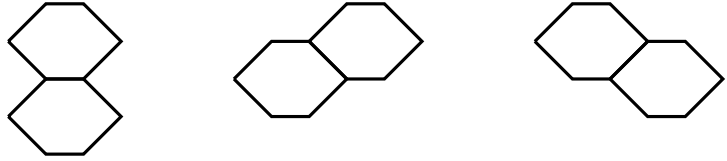
\begin{figure}[ht]
\begin{center}
\begin{tikzpicture}[scale=1]

\draw [black,very thick] (.5,.5)--(1,0)--(1.5,0)--(2,.5)--(1.5,1)--(1,1)--(.5,.5);
\draw [black,very thick] (.5,1.5)--(1,1)--(1.5,1)--(2,1.5)--(1.5,2)--(1,2)--(.5,1.5);

\draw [black,very thick] (3.5,1)--(4,0.5)--(4.5,0.5)--(5,1)--(4.5,1.5)--(4,1.5)--(3.5,1);
\draw [black,very thick] (4.5,1.5)--(5,1)--(5.5,1)--(6,1.5)--(5.5,2)--(5,2)--(4.5,1.5);

\draw [black,very thick] (7.5,1.5)--(8,1)--(8.5,1)--(9,1.5)--(8.5,2)--(8,2)--(7.5,1.5);
\draw [black,very thick] (8.5,1)--(9,0.5)--(9.5,0.5)--(10,1)--(9.5,1.5)--(9,1.5)--(8.5,1);

\end{tikzpicture}

\caption[]{The three hexagonal polyominoes with two hexagons.}

\label{Fig1}
\end{center}
\end{figure}

The number of nonequivalent hexagonal polyominoes with $n$ cells, $a(n)$, is given by Sloane's sequence
{\bf A001207} for $n \leq 35$ [NJAS]. For example, $a(1)=1$, $a(2)=3$, $a(3)=11$,
$a(4)=44$, $a(5)=186$, $a(6)=814$, $a(7)=3652$, $a(8)=16689$, $a(9)=77,359$, and $a(10)=362,671$.

A hexagonal polyominoes in which each column contains at most $k$ contiguous blocks of
cells is referred to as a $k$-board polyominoes.  Figure 2 shows examples of one and two-board polyominoes
 with $7$ hexagons.  The number of nonequivalent one-board polyominoes with $n$ hexagonal cells, $b(n)$, is given by Sloane's sequence {\bf A059716} for $n \leq 24$ [GV]. For example, $b(1)=1$, $b(2)=3$,
$b(3)=11$, $b(4)=42$, $b(5)=162$, $b(6)=626$, $b(7)=2419$, $b(8)=9346$,
$b(9)=36,106$, and $b(10)=139,483$.

\begin{figure}[ht]
\begin{center}
\begin{tikzpicture}[scale=1]

\draw [black,very thick] (-.5,1.5)--(0,1)--(.5,1)--(1,1.5)--(.5,2)--(0,2)--(-.5,1.5);
\draw [black,very thick] (.5,1)--(1,0.5)--(1.5,0.5)--(2,1)--(1.5,1.5)--(1,1.5)--(.5,1);
\draw [black,very thick] (.5,2)--(1,1.5)--(1.5,1.5)--(2,2)--(1.5,2.5)--(1,2.5)--(.5,2);
\draw [black,very thick] (1.5,.5)--(2,0)--(2.5,0)--(3,.5)--(2.5,1)--(2,1)--(1.5,.5);
\draw [black,very thick] (1.5,1.5)--(2,1)--(2.5,1)--(3,1.5)--(2.5,2)--(2,2)--(1.5,1.5);
\draw [black,very thick] (1.5,2.5)--(2,2)--(2.5,2)--(3,2.5)--(2.5,3)--(2,3)--(1.5,2.5);
\draw [black,very thick] (2.5,1)--(3,0.5)--(3.5,0.5)--(4,1)--(3.5,1.5)--(3,1.5)--(2.5,1);

\draw (2,-.7)  node {(a)};

\draw [black,very thick] (10,1.5)--(10.5,1)--(11,1)--(11.5,1.5)--(11,2)--(10.5,2)--(10,1.5);
\draw [black,very thick] (11,1)--(11.5,0.5)--(12,0.5)--(12.5,1)--(12,1.5)--(11.5,1.5)--(11,1);
\draw [black,very thick] (11,3)--(11.5,2.5)--(12,2.5)--(12.5,3)--(12,3.5)--(11.5,3.5)--(11,3);
\draw [black,very thick] (12,.5)--(12.5,0)--(13,0)--(13.5,.5)--(13,1)--(12.5,1)--(12,.5);
\draw [black,very thick] (12,1.5)--(12.5,1)--(13,1)--(13.5,1.5)--(13,2)--(12.5,2)--(12,1.5);
\draw [black,very thick] (12,2.5)--(12.5,2)--(13,2)--(13.5,2.5)--(13,3)--(12.5,3)--(12,2.5);
\draw [black,very thick] (13,1)--(13.5,0.5)--(14,0.5)--(14.5,1)--(14,1.5)--(13.5,1.5)--(13,1);

\draw (12,-.7)  node {(b)};

\end{tikzpicture}

\caption[]{A one-board hexagonal polyominoes (a) and a two-board hexagonal polyominoes (b) with $7$ hexagons.}

\label{Fig1}
\end{center}
\end{figure}
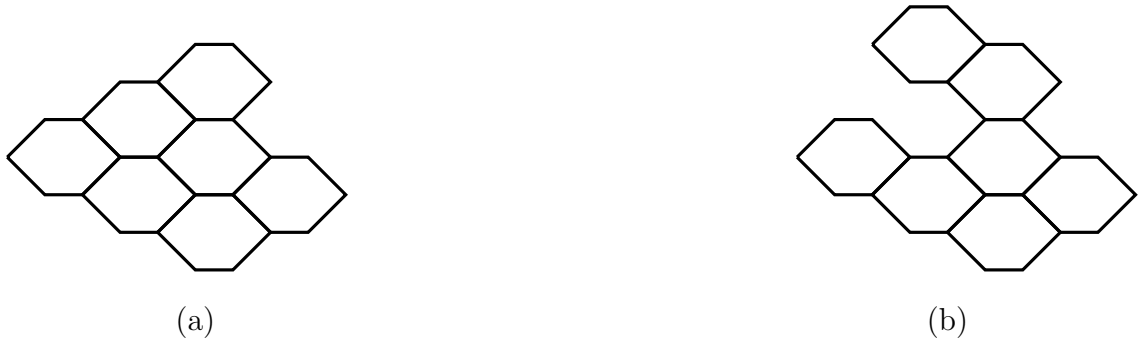

The generating function that enumerates one-board hexagonal polyominoes is computed in [K], which can also be computed using the Maple package LEGO in [Z0] by taking $p(a,b)=a+b$ (for definition of $p(a,b)$ refer to [Z0]). To the best of the authors' knowledge, $k$-board hexagonal polyominoes have not been enumerated for $k \geq 2$. In this article, in addition to computing generating functions (and series expansions) for hexagonal lattice animals that fit into a prescribed but arbitrary height, we also compute the first $12$ terms of the sequence that enumerate the number of board-pair-hexagonal polyominoes, the analog of square board-pair-polyominoes computed in [Z2]. Fisrt we recall some definitions and introduce the essential details of the Maple implementation of the transfer-matrix method [S] from [Z1].

\vskip 10pt

\section*{\normalsize 2. The Transfer-Matrix Method}
Let $G(V,E)$ be a vertex weighted directed graph consisting a finite set of vertices, $V$, and a finite set of directed edges, $E$. A path $P$ in G is a sequence $v_1,e_1,v_2,e_2,\dots,v_{m-1},e_m,v_m$
where $v_i \in V$($1 \leq i \leq m$) and $e_i \in E$($1 \leq i \leq m$) is
an edge from $v_{i-1}$ to $v_{i}$. The weight of a path $P$ is the sum of
the weights of all vertices participating in the path $P$. In this case $Wt(P)=Wt(v_1)+Wt(v_2)+\ldots +Wt(v_m)$.

For $T$ a subset of $V$, we want to compute the weight enumerator (generating function),

\begin{equation}
F(z)= \sum_{j=0}^{\infty} a_jz^j \nonumber
\end{equation}

\noindent where $a_j$ is the number of paths with weight $j$ that starts at any vertex in $V$ and ends at a vertex in $T$. For $v \in V$, let $F_v(z)$ be the generating function that enumerates all paths starting at vertex $v$ and ends at a vertex in $V$. Then, clearly,

\begin{equation}
  F_v(z)=\sum_{\substack {paths \,\,\, \textbf{P} \\  init(\textbf{P})=v, \,\, fin(\textbf{P}) \in T}}t^{wt(\textbf{P})}
\end{equation}

\noindent where $init(P)$ is the initial vertex of $P$ and $fin(P)$ is the
terminal vertex of $P$. Let $N(v)$ be the set of vertices adjacent to $v$.
Then, $\set{F_v(z)|v \in V}$ satisfies the following $\norm{V}$ (where $\norm{V}$ is the number of vertices) system of
equations in $\norm{V}$ unknowns, namely $\set{F_v(z)|v \in V}$,

\begin{equation}
F_v(z)={\bf 1}_T (v)z^{wt(v)} + z^{wt(v)}\sum_{u \in N(v)}
n(v,u) F_u(z)
\end{equation}

\noindent where $n(v,u)$ is the number of paths from $v$ to $u$ and

\begin{equation}
 \bf{1}_T(z)= \left\{
\begin{array}{cc}
1 & if  \,\, \,\,z \in T\\
0 & \,\,\,\, otherwise
\end{array} \right.
\end{equation}

Since each $F_u(z)$ exists as a power series in $z$, our system has a solution. Once we solve for this sytem, the answer to our original problem is then

\begin{equation}
F(z)=\sum_{v \in T} F_v(z) \nonumber
\end{equation}

As an example consider the following vertex-weighted directed graph on four vertices with weights as shown next to the vertices.

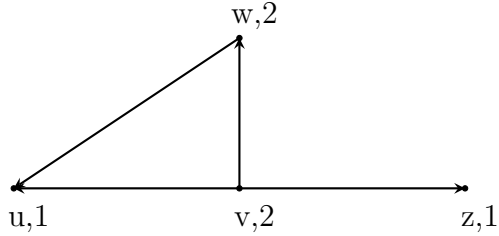
\begin{figure}[ht]
\begin{center}

\begin{tikzpicture}[>=stealth][scale=1]
\filldraw [black] (0,0) circle (1pt);
\filldraw (3,0) circle (1pt);
\filldraw (3,2) circle (1pt);
\filldraw (6,0) circle (1pt);
\draw (0.2,-0.4)  node {u,1};
\draw (3.2,-0.4)  node {v,2};
\draw (6.2,-0.4)  node {z,1};
\draw (3.2,2.3)  node {w,2};
\draw [<-, thick] (0,0)-- (3,2);
\draw [->, thick] (3,0)--(3,2);
\draw [<-, thick] (0,0)-- (3,0);
\draw [->, thick] (3,0)--(6,0);
\draw (0,0)-- (3,0) -- (6,0);
\end{tikzpicture}

\caption{Example of vertex weighted directed graph on four vertices}
\label{Fig 1}

\end{center}
\end{figure}

In this example, the generating function that enumerates all paths according to their weight is found by first finding the four generating functions $\{F_u(z),F_v(z),F_w(z),F_t(z)\}$ for all paths starting at vertex $u, v, w,$ and $t$ respectively. These four generating functions are related by $$ f_t(z)=z,f_w(z)=f_v(z)=zf_u(z)=\frac{z(1+z+z^2+z^3+z^4)}{1-z^5}\,.$$
Solving these system and adding together we get the required generating function that enumerates all paths according to their weight as

$$
F(z)=\frac{2z+2z^2+z^3+3z^4+3z^5+z^6}{1-z^5} \, .
$$

For the hexagonal polyominoes case, we follow the structure of [Z1] and automate the problem and use
Maple do the hard part. First, we recall the following definitions from
[Z1]:

{\bf Definition 2.1} A {\it Combinatorial Markov Process} is a six tuple {\mbox
(V, E, init, fin, Start, Finish)}, where $V$ is a finite set of vertices, $E$
is the set of directed edges, and $init, fin: E \To V$ are functions that
assign the initial and terminal vertex to an edge, respectively.

{\bf Definition 2.2} A vertex {\it Weighted Combinatorial Markov Process}
is a seven tuple {\mbox (V, E, init, fin, Start, Finish, wt)}, where {\mbox (V,
E, init, fin, Start, Finish)} is a Combinatorial Markov Process and $wt: V(E)
\To \posint$ is a function that assigns weights to vertices.

In Maple, we represent a Combinatorial Markov Process with vertices
$V=\set{1,2,\dots,n}$ and no multiple edges as a four tuple {\tt
[S,T,ListOfOutGoingNeighbors,ListOfWeights]} where {\tt S},{\tt T} $\subseteq
V$, and {\tt ListOfOutGoingNeighbors} and {\tt ListOfWeights} are lists of
length $n$ whose $i^{th}$ element is the set of vertices adjacent to vertex $i$
and the weight of vertex $i$, respectively.

Similarly, the Combinatorial Markov Process for a multiple edge graph is a four
tuple {\tt [S,T,ListOfOutGoingNeibors,ListOfWeights]} where {\tt S}, {\tt T},
and {\tt ListOfWeights} are as previously defined and {\tt
ListOfOutGoingNeighbors} is a list of length $n$ whose $i^{th}$ element is a
multi-set $j_1^{m_1} j_2^{m_2} \dots j_k^{m_k}$ represented in the form
$\set{[j_1,m_1],[j_2,m_2],\dots,[j_k,m_k]}$, which means that from vertex $i$
there are $m_1$ edges to vertex $j_1$, $m_2$ edges to vertex $j_2$, and so on.

If we successfully model hexagonal polyominoes as a weighted Combinatorial Markov Process as it is done for the squre lattice animals in [Z1] then we can employ the Maple package {\tt MARKOV} given in [Z1] to automatically find the generating function and series expansion of the generating function. The next section describes how weighted Combinatorial Process for hexagonal polyominoes is created.

\vskip 10pt

\section*{\normalsize 3. Maple Representation of hexagonal polyominoes}

We could develop a grammar that describes height-restricted hexanimals from first principles
and apply the same methods as Zeilberger [Z1] to encode the grammar as a Combinatorial
Markov Process; however, this would result in rewriting much of the code that
was already created for the square lattice animals. Our approach here is to define a bijection between hexagonal lattice animals and a special class of square lattice animals and adopt the Maple codes ANIMALS and FreeANIMALS in [Z1]. In the remaining sections, we describe this processes.

We overlay a hexagonal polyominoe on the square lattice so that its leftmost vertex (vertices of the square contained inside the hexagonal) lie on the line $x=0$ and its bottom most cell lie on the line $y=0$ as shown in Fig.4(a). We then construct a polyominoes comprised of the square cells which have diagonally opposing corners falling on the same hexagonal cell as shown in Fig. 4(b).

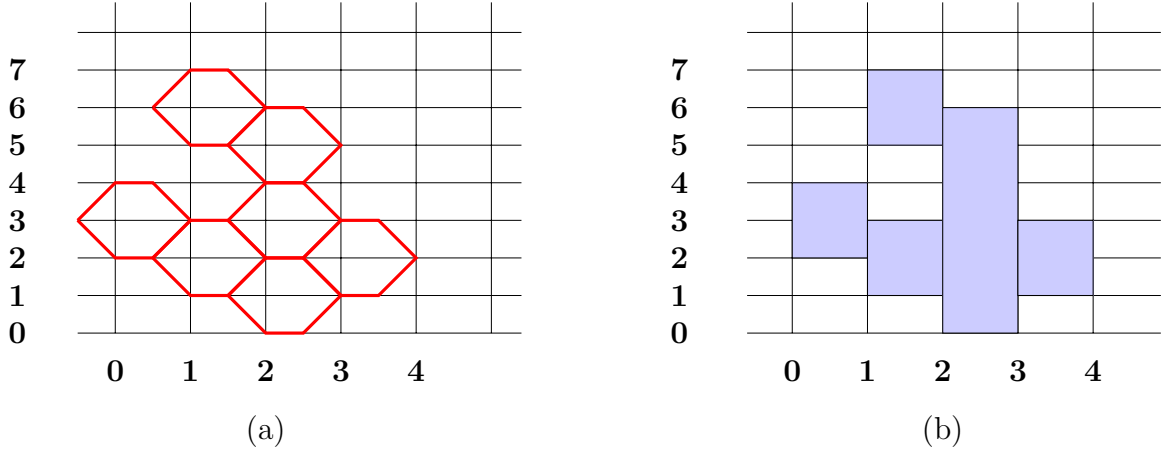
\begin{figure}[ht]
\begin{center}
\begin{tikzpicture}[scale=1]

\draw [xstep=1cm,ystep=.5cm, very thin] (-.5,0) grid (5.4,4.4);

\draw (0,-.5)  node {\textbf{0}};
\draw (1,-.5)  node {\textbf{1}};
\draw (2,-.5)  node {\textbf{2}};
\draw (3,-.5)  node {\textbf{3}};
\draw (4,-.5)  node {\textbf{4}};

\draw (-1.3,0)  node {\textbf{0}};
\draw (-1.3,0.5)  node {\textbf{1}};
\draw (-1.3,1)  node {\textbf{2}};
\draw (-1.3,1.5)  node {\textbf{3}};
\draw (-1.3,2)  node {\textbf{4}};
\draw (-1.3,2.5)  node {\textbf{5}};
\draw (-1.3,3)  node {\textbf{6}};
\draw (-1.3,3.5)  node {\textbf{7}};

\draw [red,very thick] (-.5,1.5)--(0,1)--(.5,1)--(1,1.5)--(.5,2)--(0,2)--(-.5,1.5);
\draw [red,very thick] (.5,1)--(1,0.5)--(1.5,0.5)--(2,1)--(1.5,1.5)--(1,1.5)--(.5,1);
\draw [red,very thick] (.5,3)--(1,2.5)--(1.5,2.5)--(2,3)--(1.5,3.5)--(1,3.5)--(.5,3);
\draw [red,very thick] (1.5,.5)--(2,0)--(2.5,0)--(3,.5)--(2.5,1)--(2,1)--(1.5,.5);
\draw [red,very thick] (1.5,1.5)--(2,1)--(2.5,1)--(3,1.5)--(2.5,2)--(2,2)--(1.5,1.5);
\draw [red,very thick] (1.5,2.5)--(2,2)--(2.5,2)--(3,2.5)--(2.5,3)--(2,3)--(1.5,2.5);
\draw [red,very thick] (2.5,1)--(3,0.5)--(3.5,0.5)--(4,1)--(3.5,1.5)--(3,1.5)--(2.5,1);

\draw (2,-1.3)  node {(a)};


\draw [xstep=1cm,ystep=.5cm, very thin] (8.4,0) grid (13.9,4.4);

\draw (9,-.5)  node {\textbf{0}};
\draw (10,-.5)  node {\textbf{1}};
\draw (11,-.5)  node {\textbf{2}};
\draw (12,-.5)  node {\textbf{3}};
\draw (13,-.5)  node {\textbf{4}};

\draw (7.5,0)  node {\textbf{0}};
\draw (7.5,0.5)  node {\textbf{1}};
\draw (7.5,1)  node {\textbf{2}};
\draw (7.5,1.5)  node {\textbf{3}};
\draw (7.5,2)  node {\textbf{4}};
\draw (7.5,2.5)  node {\textbf{5}};
\draw (7.5,3)  node {\textbf{6}};
\draw (7.5,3.5)  node {\textbf{7}};

\filldraw[fill=blue!20] (9,1) rectangle (10,2);
\filldraw[fill=blue!20] (10,3.5) rectangle (11,2.5);
\filldraw[fill=blue!20] (10,1.5) rectangle (11,.5);
\filldraw[fill=blue!20] (11,0) rectangle (12,3);
\filldraw[fill=blue!20] (12,0.5) rectangle (13,1.5);

\draw (11,-1.3)  node {(b)};

\end{tikzpicture}

\caption[]{Mapping between a hexagonal lattice animals and a square lattice animals}

\label{Fig1}
\end{center}
\end{figure}

As described in detail by Zeilberger [Z1], we represent a polyomino by the
set of coordinates of the bottomleft corner of each of its cells.  We further
encode a polyomino as a word in the alphabet consisting of the non-empty
subsets of $\set{0,1,...,k}$ where $k$ is the maximum $y$ coordinate in the set
representation of the polyomino.  For example, the animal in Fig. 4
(b) is represented by the set
$$
\set{(0,2),(0,3),(1,1),(1,2),(1,5),(1,6),(2,0),(2,1),(2,2),(2,3),(2,4),(2,5),(3,1),(3,2)}
$$
\noindent and encoded as the word (in interval notation)
$$
\set{[2,3]\},\{[1,2],[5,6]\},\{[0,5]\},\{[1,2]}.
$$

There are two important properties of the polyomino words that encode
hexanimals.  First, every interval has even size (i.e., it represents a
contiguous block of an even number of cells); this is obvious from the mapping
which associates each hexagonal cell with a pair of square cells.  Second,
within a letter, the starting position of every interval has the same parity
(odd or even), and adjacent letters have different parity; this is again
obvious from the mapping of each hexagonal cell to a pair of square cells and
the vertical shift of adjacent hexagonal columns.  Note, however, that the
parity of the leftmost column may be either odd or even depending on the
hexanimal.

We will call the class of square polyominoes to which hexagonal polyominoes are mapped {\it
parity polyominoes}.  It is clear that the mapping between hexagonal animals and
parity polyominoes is unique and reversible, and thus defines a bijection.  As
a result of this bijection, we can count hexagonal animals by counting parity
polyominoes. Zeilberger's code counts general polyominoes (a much larger class
than parity polyominoes), so in order to use it to count parity polyominoes, we
need to restrict the alphabet.  In the following sections of this paper, we
describe how to modify Zeilberger's code to count parity polyominoes while
minimizing the amount of changes required; we assume that the reader is
familiar with [Z1].

\vskip 10pt

\section*{\normalsize 4. Counting Globally Skinny Hexagonal Animals}
Zeilberger's {\tt ANIMALS} package uses the {\it transfer-matrix method} [S] to count polyominoes whose set representation has $y$ coordinates in the range $0
\leq y \leq k$ for some arbitrary $k$. Before the transfer-matrix method can be
applied, a grammar needs to be constructed to describe how to combine letters
into words that encode polyominoes.  A key routine in the generation of this
grammar is {\tt PreLeftLetters(a,b)} which returns the set of all possible
subsets of the integers in the range $[a,b]$ written in interval notation. For
example, {\tt PreLeftLetters(0,2)} returns the following subsets:
$$
\set{\set{}, \set{[0,0]}, \set{[0,1]}, \set{[1,1]}, \set{[0,2]},
\set{[1,2]}, \set{[2,2]}, \set{[0,0], [2,2]}}
$$

For parity polyominoes, we are only interested in those subsets whose intervals
have even size and the same parity for the starting position.  Rather than
rewriting \linebreak {\tt PreLeftLetters} from scratch, we can use a trick in
Maple to subclass {\tt PreLeftLetters} and alter its parameters without having
to modify the original {\tt PreLeftLetters} code:

\begin{verbatim}
read ANIMALS:

origPreLeftLetters := subs(PreLeftLetters=origPreLeftLetters,
                           eval(PreLeftLetters)):

PreLeftLetters := proc(a, b)
    local halfS, S, i, PreLet,  PreLet1:

    halfS := origPreLeftLetters(0, floor((b-a+1)/2)-1) minus {{}}:
    S := {{}}:
    for i from 1 to nops(halfS) do
        PreLet  := {seq([a+2*halfS[i][j][1], a+2*halfS[i][j][2]+1],
                    j=1..nops(halfS[i]))}:
        PreLet1 := {seq([PreLet[j][1]+1, PreLet[j][2]+1],
                    j=1..nops(PreLet))}:
        S := S union {PreLet}:
        if max(seq(PreLet1[j][2], j=1..nops(PreLet1))) <= b then
            S := S union {PreLet1}:
        fi:
    od:
    S:
end:
\end{verbatim}

The first line loads Zeilberger's original {\tt ANIMALS} package. The second
line makes a copy of {\tt PreLeftLetters}; the {\tt subs} function is designed
to modify {\tt PreLeftLetters} so that any recursive calls are made to the copy
rather than the original. The remaining lines redefine {\tt PreLeftLetters} so
that it produces the type of intervals that are required by parity polyominoes.
The new {\tt PreLeftLetters} uses the original {\tt PreLeftLetters} to generate
a set of subsets (in interval notation) for half the required range (stored in
{\tt halfS}). The {\tt for} loop goes through each interval and doubles its
size; this results in intervals that have even size and even parity (stored in
{\tt PreLet}).  To create intervals with even size and odd parity, each
interval in {\tt PreLet} is shifted down by one, and the result is stored in
{\tt PreLet2}.  The final condition checks that the shifted intervals do not go
out of the range $[a,b]$.

For example, {\tt PreLeftLetters(0,5)} would call {\tt origPreLeftLetters(0,2)}
and get the subsets shown above.  The {\tt for} loop would then modify these
subsets to generate the following result which represents the parity polyomino
preletters in the range $[0,5]$:

\begin{small}
$$
\set{\set{}, \set{[0,1]}, \set{[1,2]}, \set{[0,3]}, \set{[1,4]}, \set{[2,3]},
\set{[3,4]}, \set{[0,5]}, \set{[2,5]}, \set{[4,5]}, \set{[0,1],[4,5]}}
$$
\end{small}

Once the preletters are generated, the grammar is extended by determining which
of the preletters are valid extensions of existing polyomino words; the {\tt
ANIMALS} routine that performs this check is {\tt PreLetToLet}.  For parity
polyominoes, not only must the preletter satisfy the same connectivity
constraints as for general polyominoes, but the parity of the preletter's
interval starting positions must be different than the letter it follows.  As
we did with {\tt PreLeftLetters}, a simple extension to the original code is
that is required:

\begin{verbatim}
origPreLetToLet := subs(PreLetToLet=origPreLetToLet,
                        eval(PreLetToLet)):

PreLetToLet := proc(Let, PreLet)
        if PreLet <> {} and Let[1][1][1] mod 2 = PreLet[1][1] mod 2 then
                0:
        else
                origPreLetToLet(Let, PreLet):
        fi:
end:
\end{verbatim}

Since we know that every interval in a parity polyomino letter has the same
starting position parity, it suffices to check whether or not the starting
positions of the first intervals have differing parities; if that is the case,
{\tt PreLetToLet} calls the original {\tt ANIMALS} routine to ensure that the
connectivity constraints are satisfied.

In order to weight the transfer-matrix vertices so that polyominoes with a
specific number of cells can be enumerated, the {\tt ANIMALS} routines call
{\tt Weight(Let)} to determine the number of cells in the letter {\tt Let}.
For enumerating hexanimals, we want to weight the transfer-matrix vertices with
the number of hexagonal cells.  Since each hexagonal cell is represented by two
polyomino cells, a letter's weight in hexagonal cells is half its weight in
polyomino cells:

\begin{verbatim}
origWeight := subs(Weight=origWeight,
                   eval(Weight)):

Weight := proc(Let)
        floor(origWeight(Let)/2):
end:
\end{verbatim}

One final change to {\tt ANIMALS} is needed in order to enumerate hexanimals;
the routine {\tt Khaya(L)} that computes the number of polyominoes with $\leq
L$ cells relies on symmetry that does not exist in the hexagonal lattice.
Since a hexanimal with $L$ cells must map to a parity polyomino whose set
representation has a maximum $y$ coordinate of $2L$, we can implement {\tt
Khaya(L)} by a call to {\tt GFseries}:

\begin{verbatim}
Khaya := proc(L)
        [GFseries(2*L, L)]:
end:
\end{verbatim}

\vskip 10pt

\section*{\normalsize 5. A User's Manual for HexANIMALS}
The modifications to {\tt ANIMALS} that are described in the previous section
are contained in the package {\tt HexANIMALS}. To run {\tt HexANIMALS},
download it and {\tt ANIMALS} from the first authors' web site to a local directory
and start Maple.  Once in Maple, type: {\tt read HexANIMALS} and follow the
on-line help.

Excluding {\tt Khaya}, all the generating function and series expansion
routines in {\tt ANIMALS} remain unchanged. {\tt GF(n,s)} computes the
generating function for hexanimals embedded in the square lattice of height $n$
(i.e., the height of the square lattice is comprised of $n$ cells).  {\tt
GFseries(n,L)} computes the list of length $L$ whose $k^{th}$ term is the
number of hexanimals with $k$ cells embedded in the square lattice of height
$n$.  {\tt Khaya(L)} computes the list of length $L$ whose $k^{th}$ term is the
number of hexanimals with $k$ cells.  {\tt Gf} and {\tt Gfseries} are the
analogs of {\tt GF} and {\tt GFseries} for hexanimals whose polyomino
embeddings have height {\it exactly} $n$.

As an example of {\tt HexANIMALS} usage, the call {\tt GF(1,s)} returns $0$
since no hexanimal can be embedded in a square lattice of height $1$. {\tt
GF(2,s)} returns $s$ since there is exactly one hexanimal (the single cell)
that can be embedded in a square lattice of height $2$. {\tt GFseries(3,5)}
returns $[1,2,2,2,2]$ since the hexanimals in the square lattice of height $3$
with more than one hexagonal cell are chains, and each chain has two
orientations depending on the parity of the starting cell. In Table
\ref{tab:global} we compute the generating function for hexanimals embedded in
a square lattice of height $n \leq 7$ along with the first few terms in the
expansion by hexagonal cell count.  Once $n\geq 13$, it takes too long to
compute {\tt GF(n,s)} exactly, but one can go much further with {\tt
GFseries(n,L)}.

\begin{table}[!h]
\begin{center}
\begin{tabular}{l|l|l}
\hline
$n$ & $F(z)$ & $a_i$\\
\hline
$1$ & $0$ & $0,0,0,\dots$\\
\hline
$2$ & $z$ & $1,0,0,\dots$\\
\hline
$3$ & $\frac{z(1+z)}{1-z}$ & $1,2,2,\dots$\\
\hline
$4$ & $\frac{z(1+z)}{(z^2+z-1)(z-1)}$ & $1,3,6,11,19,32,53,87$\\
\hline
$5$ & $\frac{z(z^3+2z^2+z+1}{(z^3+2z^2+z-1)(z^2+z-1)}$ & $1,3,10,25,61,142,323,723$\\
\hline $6$ & $\frac
{z(1+z+2z^2+7z^5+4z^3+8z^4-3z^6-4z^7+3z^8+z^10+4z^9)}{(3z^6+2z^7-1+2z^2+z+5z^3)(z^3+2z^2+z-1)}$ & $1,3,11,37,111,320,896$\\
\hline
\end{tabular}
\caption[]{Generating functions for hexanimals polyominoes embedded in the square
lattice of height $n$ and their series expansion by hexagonal cell count.}
\label{tab:global}
\end{center}
\end{table}

\vskip 10pt

\section*{\normalsize 6. Counting Locally Skinny Hexagonal Animals}
If we relax the condition that the entire hexagonal polyominoes must fit within a square
lattice of fixed height, and only require that each column have height less
than some arbitrary constant, then we get what Zeilberger calls locally skinny
hexagonal polyominoes (as opposed to globally skinny hexagonal polyominoes).

The code for computing the generating function and series expansions for
locally skinny polyominoes is contained in Zeilberger's {\tt FreeANIMALS}
package.  Since the changes to {\tt FreeANIMALS} that are required in order to
enumerate hexagonal lattice animals are virtually the same as those made to {\tt ANIMALS}, we
omit the details that were discussed in the previous sections and only discuss
the major differences.

{\tt FreeANIMALS} creates a Combinatorial Markov Process with multiple edges
between letters.  Each edge represents one of the possible vertical offsets of
the letters.  In contrast, {\tt ANIMALS} creates a Combinatorial Markov Process
with only one edge between letters; thus, each letter represents a particular
arrangement of cells with an intrinsic offset.  As a result, {\tt
PreLeftLetters} in {\tt HexANIMALS} had to create the two possible parity
offsets for each letter whereas {\tt PreLeftLetters} in {\tt HexaFreeANIMALS}
needs to only create a single instance of the letter and the code that
generates the Combinatorial Markov Process will create its own instances of the
shifted letter.  The resulting {\tt PreLeftLetters} routine is similar to the
{\tt PreLeftLetters} routine in {\tt HexANIMALS} except that the parity
shifting code is omitted:

\begin{verbatim}
read FreeANIMALS:

origPreLeftLetters := subs(PreLeftLetters=origPreLeftLetters,
                           eval(PreLeftLetters)):

PreLeftLetters := proc(a, b)
    local halfS, S, i, PreLet:

    halfS := origPreLeftLetters(0, floor((b-a+1)/2)-1) minus {{}}:
    S := {{}}:
    for i from 1 to nops(halfS) do
        PreLet := {seq([a+2*halfS[i][j][1],
                        a+2*halfS[i][j][2]+1], j=1..nops(halfS[i]))}:
        S := S union {PreLet}:
    od:
    S:
end:
\end{verbatim}

In addition, {\tt FreeANIMALS} contains the routine {\tt
PreLeftLettersk(a,b,k)} which computes all possible subsets of the integers in
the range $[a,b]$ that are represented by exactly $k$ intervals.  {\tt
PreLeftLettersk} is used to compute $k$-board polyominoes.  For example, {\tt
PreLeftLettersk(0,4,2)} returns the following subsets:

$$
\set{\set{[0,0],[2,2]},\set{[0,0],[2,3]},\set{[0,0],[3,3]},\set{[0,1],[3,3]}}
$$

The modifications to convert {\tt PreLeftLettersk} for use with {\tt
HexaFreeANIMALS} are analogous to the changes made to {\tt PreLeftLetters} as
shown in the following code:

\begin{verbatim}
PreLeftLettersk := proc(a, b, k)
    local halfS, S, i, PreLet:

    halfS := origPreLeftLettersk(0, floor((b-a+1)/2)-1, k) minus {{}}:
    S := {{}}:
    for i from 1 to nops(halfS) do
        PreLet := {seq([a+2*halfS[i][j][1],
                   a+2*halfS[i][j][2]+1], j=1..nops(halfS[i]))}:
        S := S union {PreLet}:
    od:
    S:
end:
\end{verbatim}

The only other significant difference between {\tt HexANIMALS} and {\tt
HexFreeANIMALS} is in the specification of height restrictions.  When
considering globally skinny hexanimals, the height restrictions were defined in
terms of the square lattice in which the hexanimals were embedded; this was
because the shifting of adjacent columns in hexanimals makes the definition of
the height of a hexanimal in terms of hexagonal cells difficult.  In the case
of locally skinny hexanimals, the height restrictions apply to individual
columns and are well-defined in terms of hexagonal cells.  Since the underlying
algorithms of {\tt HexaFreeANIMALS} operate on parity polyominoes, the height
restrictions for hexanimal columns are simply doubled when applied to parity
polyominoes as shown in the following code:

\begin{verbatim}
origgf := subs(gf=origgf, eval(gf)):

gf := proc(n, s)
    origgf(2*n, s):
end:

origgfSeries := subs(gfSeries=origgfSeries, eval(gfSeries)):

gfSeries := proc(n, L)
    origgfSeries(2*n, L):
end:

origgfList := subs(gfList=origgfList, eval(gfList)):

gfList := proc(resh, s)
    origgfList([seq(2*resh[i], i=1..nops(resh))], s):
end:

origgfSeriesList := subs(gfSeriesList=origgfSeriesList,
                         eval(gfSeriesList)):

gfSeriesList := proc(resh, L)
    origgfSeriesList([seq(2*resh[i], i=1..nops(resh))], L):
end:
\end{verbatim}

\vskip 10pt

\section*{\normalsize 7. A User's Manual for the HexaFreeANIMALS}
The modifications to {\tt FreeANIMALS} are contained in the package {\tt
HexaFreeANIMALS}. \linebreak {\tt gf(n,s)} computes the generating function for
hexagonal polyominoes whose columns span $\leq n$ hexagonal cells.  {\tt gfSeries(n,L)}
computes the list of length $L$ whose $k^{th}$ term is the number of hexagonal polyominoes
with $k$ cells whose columns span $\leq n$ hexagonal cells.  {\tt
gfList(List,s)} computes the generating function for polyominoes whose $k$-board
columns span $\leq List[k]$ hexagonal cells. {\tt gfSeriesList(List,L)}
computes the series expansion of \linebreak {\tt gfList(List,s)} up to $L$
terms.

As an example of the usage of {\tt HexaFreeANIMALS}, the call {\tt
gfList([7,5],s)} would compute the generating function for polyominoes whose
$1$-board columns span $\leq 7$ hexagonal cells and whose $2$-board columns
span $\leq 5$ hexagonal cells.  {\tt gfSeriesList([24],24)} computes Sloane's
sequence {\bf A0059716} [NJAS] of $1$-board polyominoes.

The call {\tt gfSeriesList([12,12],12)} computes the first $12$ terms of the
previously unknown sequence

$$
1,3,11,44,186,814,3648,16611,76437,354112,1647344,7682237
$$

\noindent that enumerates $2$-board hexanimals with a given hexagonal cell
count, the analog of Sloane's sequence {\bf A001170} for board-pair-pile
polyominoes.

\vskip 5pt

\section*{\normalsize 8. Conclusion}
We have described the Maple packages {\tt HexANIMALS} and {\tt
HexaFreeANIMALS}. {\tt HexANIMALS} enumerates globally skinny hexagonal polyominoes in
which each polyominoes is embedded in a horizontal strip of prescribed height.
{\tt HexaFreeANIMALS} enumerates locally skinny polyominoes in which the total
height of each hexagonal polyominoes is unbounded but whose columns are composed of a
bounded number of hexagonal cells.  In addition, we also demonstrated used
{\tt HexaFreeANIMALS} to compute the first few terms of the number of board-pair-hexagonal animals.

\vskip 5pt

\section*{\normalsize References}

\noindent [K] David A. Klarner, {\it Cell Growth Problems}, Canad. J. Math.
{\bf 19} (1967), pp. 851-863.

\noindent [GV] Markous Voge and Anthony J. Guttmann, {\it On the number of hexagonal polyominoes}, Theoret. Comp. Sci. 307 (2003), 433-53.

\noindent [NJAS] N.J. Sloane and S. Plouffe, The Encyclopedia of Integer Sequences, Academic Press, San Diego, 1995.

\noindent [S] Richard P. Stanley, "Enumerative Combinatorics", Vol. 1,
Wadsworth, Monterey, California, 1986.

\noindent [Z0] Doron Zeilberger, {\it Automated Counting of LEGO Towers}, J.
Difference Eq. {\bf 5} (1999), pp. 355-377.

\noindent [Z1] Doron Zeilberger, {\it Symbol-Crunching with the Transfer-Matrix
Method in Order to Count Skinny Physical Creatures}, INTEGERS ({\tt
http://www.integers-ejcnt.org}) {\bf 0} (2000), A9.

\noindent [Z2] Doron Zeilberger, {\it The Umbral Transfer-Matrix Method III,
Counting Animals}, New York J. of Mathematics {\bf 7} (2001), pp. 223-231.

\noindent
\end{document}